\documentclass[11pt]{article}
\usepackage[margin=1.0in]{geometry}
\usepackage{amsmath,amsthm,amsfonts,amssymb,amscd}
\usepackage{cite}
\usepackage{url}
\usepackage{cite}
\usepackage{fullpage}
\theoremstyle{plain}
\newtheorem{theorem}{Theorem}[section]
\newtheorem{lemma}[theorem]{Lemma}
\newtheorem{corollary}[theorem]{Corollary}
\newtheorem{proposition}[theorem]{Proposition}
\newtheorem{definition}{Definition}[section]
\newtheorem{remark}{Remark}[section]
\newtheorem{example}{Example}[section]

\usepackage{cite}
\theoremstyle{definition}

\newcommand{\la}{\langle}
\newcommand{\ra}{\rangle}

\newcommand{\nexto}{\kern -0.54em}

\newcommand{\dZ}{{\cal Z \kern -0.7em Z}}
\newcommand{\dC}{{\rm\hbox{C \kern-0.8em\raise0.2ex\hbox{\vrule height5.4pt
width0.7pt}}}}
\newcommand{\dQ}{{\rm\hbox{Q \kern-0.85em\raise0.25ex\hbox{\vrule height5.4pt
width0.7pt}}}}

\usepackage{bbm}
\usepackage{epsfig}
\usepackage{latexsym}

\newcommand{\NN}{\mathbb{N}}

\newcommand{\RR}{\mathbbm{R}}

\newcommand{\dsty}{\displaystyle}

\begin{document}
\title{Subgradient algorithms for solving variable inequalities
\thanks{This work was
partially supported by CNPq grants 303492/2013-9, 474160/2013-0 and
202677/2013-3 and by project CAPES-MES-CUBA 226/2012.} }

\author{J.Y. Bello Cruz\footnote{(Corresponding author) Instituto de Matem\'atica e Estat\'istica, Universidade Federal de Goi\'as. CEP 74001-970 GO,
Goi\^ania, Brazil. E-mail: yunier@impa.br}\and G. Bouza Allende
\footnote{Departamento de Matem\'atica Aplicada, Facultad de
Matem\'atica y Computaci\'on, Universidad de La Habana. CEP 10400
Habana, Cuba. E-mail: gema@matcom.uh} \and L.R. Lucambio
P\'erez\footnote{Instituto de Matem\'atica e Estat\'istica,
Universidade Federal de Goi\'as. CEP 74001-970 GO, Goi\^ania,
Brazil. E-mail: lrlp@mat.ufg.br}}

\maketitle

\begin{abstract} In this paper we consider the variable inequalities problem, that is, to find a solution of
the inclusion given by the sum of a function and a point-to-cone
application. This problem can be seen as a generalization of the
classical  inequalities problem taking a variable order
structure. Exploiting this relation, we propose two
variants of the subgradient algorithm for solving the variable inequalities
model. The convergence analysis is given under
convex-like conditions, which, when the point-to-cone application is
constant, contains the old subgradient schemes.
\end{abstract}
\medskip


\medskip


\section{Introduction}

 We  consider the  inclusion problem of finding $x\in C$ such that \begin{equation}\label{0inT}0\in T(x),\end{equation}
where  $T\colon\mathbb{R}^n\rightrightarrows\mathbb{R}^m$ is a point-to-set operator and $C$ is a nonempty and closed subset of
$\RR^n$. Inclusions has been studied in many
works due its applications; see, for instance, \cite{inclusion1, inclusion2,rocka1}. However, we will focus in the case in which
  $T(x)=F(x)+K(F(x))$, where
$F:\mathbb{R}^n\to\mathbb{R}^m$ and
$K\colon\RR^m\rightrightarrows\RR^m$ is a point-to-set application
such that $K(y)$ is a closed pointed convex cone for all $y\in
\RR^m$.  Then, we are lead to the model:
\begin{equation}\label{oooo}
\text{ find a point $x\in C$ fulfilling
that }0\in F(x)+K(F(x)).
\end{equation}
If  $K$ is a constant application, problem \eqref{oooo} is equivalent to compute $x\in C$ such that
\begin{equation}\label{eq11}0 \in F(x) +K.
\end{equation}
This model is known as the  $K$-inequalities problem because, using the partial order   defined in   $\RR^m$ by $K$ as  $$\hat{y} \preceq_K y\quad \text{ if and only if }\quad y-\hat{y}\in K,$$
 problem \eqref{eq11} is equivalent to:
\begin{equation}\label{pppp}
\text{find }x\in C\text{ such that } F(x)\preceq_K 0.
\end{equation}
Model
\eqref{oooo} can be interpreted as a system of  variable inequalities. Indeed,  consider the variable order  given by
$$z\preceq_{K(z)}y\quad \text{ if and only if }\quad y-z\in K(z);$$ see \cite{ap11,Gabriele-book} for more details. Then,
problem
\eqref{oooo} is equivalent to:
\begin{equation}\label{prob-1}\text{find }x\in C\text{ such that }
F(x) \preceq_{K(F(x))} 0.
\end{equation}
That is why, from now on, this problem will be called the variable inequalities problem. The solution set of this problem will be denoted by $S_*$.

Note that if $K$ is a constant application, problem \eqref{prob-1}
leads  to model \eqref{pppp}, which has been already studied in
\cite{rob1,rob,yun-luis, luis-yun}. Moreover, if $K$ is the Pareto
cone, \emph{i.e.}, $K=\RR_{+}^m$, it is equivalent to the convex
feasibility problem, which has been well-studied in \cite{baus-borw}
and has many applications in optimization theory, approximation
theory, image reconstruction and so on; see, for instance,
\cite{poljak, neumann, censor}. The variable case is not only a
generalization of problem  \eqref{pppp}. Variable order optimization
models appear in portfolio  and medicine applications, as recently reported
in \cite{bao-mord-soub1,bao-mord-soub2,Gabriele-book}.

The  algorithms for solving problem \eqref{pppp} mainly converge
under convexity of $F$. We generalize this concept to the variable
order case as follows
\begin{equation}\label{convexidade}
\alpha F(x)+(1-\alpha)F(\hat{x})-F(\alpha x+(1-\alpha)\hat{x})\in K(F(\alpha x+(1-\alpha)\hat{x})).
\end{equation}

We want to point out that relation \eqref{convexidade} generalizes the previously defined  convexity concept to the  case in which the point-to-cone application,
$K$, is identically constant.
As in this case, if $F$ is a $K$-convex function and $C$ is a convex set,    model
 \eqref{prob-1} is also called a
$K$-convex inequalities problem.

In this paper we propose a subgradient approach for solving problem
\eqref{prob-1}, which combines a subgradient iteration with a simple
projection step, onto the intersection of $C$ with suitable
halfspaces containing the solution set $S_*$. The proposed
conceptual algorithm has two variants called  Algorithm $R$ and
Algorithm $S$. The first one is based on Robinson's subgradient
algorithm given in \cite{rob} for solving problem \eqref{pppp}. The
$S$ variant corresponds to a special modification of the subgradient
algorithms proposed in \cite{sc-1} for the scalar problem ($m=1$ and
$K=\RR_+$) and in \cite{yun-luis} for solving problem \eqref{pppp}.
The main difference between the proposed variants lies in how the
projection step is done. For the convergence of the variants, we
assume that the  set $S_*$ is nonempty and that the function $F$ is
$K$-convex with respect to the defined variable order extending the
previous schemes.

The paper is organized as follows. In the next section, we outline
the main definitions and preliminary results. In Section \ref{sec-3}
some analytical results and comparisons for $K$-convex functions are
established. Section \ref{sec-4} is devoted to the presentation of
the algorithms and their convergence is shown in Section
\ref{sec-5}. Finally, some comments and remarks are presented in
Section \ref{sec-6}.

\section{Preliminaries}\label{preliminar}
In this section, we present some definitions and results, which are
needed in the convergence analysis. We begin with some classical
notations.

The inner product in $\RR^n$ is denoted by $\la\cdot,\cdot\ra$, the
norm, induced by this inner product, by $\|\cdot\|$ and $B[x,\rho]$
is the closed ball centered at $x\in\RR^n$ with radio $\rho$,
\emph{i.e.}, $B[x,\rho]:=\{y\in\RR^n\colon\|y-x\|\leq\rho\}$. A set
valued application  $K:\RR^m\rightrightarrows \RR^m$ is closed if
and only if $gr(K):=\{(x,y)\in \RR^m\times\RR^m: y\in K(x)\}$ is a
closed set. Given the cone $\mathcal{K}$, the dual cone of
$\mathcal{K}$, denoted $\mathcal{K}^*$, is
$\mathcal{K}^*:=\left\{z\in \RR^m\colon \langle z\,,\, y\rangle\geq
0, \; \forall y\in \mathcal{K} \right\}.$

The set $C$ will be a closed and convex subset of $\RR^n$. For an
element $x\in \RR^n$, we define the orthogonal projection of $x$
onto $C$,  $P_C(x)$, as the unique point in $C$, such that
$\|P_C(x)-y\|\le\|x-y\|$ for all $y\in C$. In the following we consider a well
known fact on orthogonal projections.

\begin{proposition}\label{popiedades_projeccion}
Let $C$  be a nonempty, closed and convex set in $\RR^n$. For all
$x\in \RR^n$ and all $z\in C$, the following property holds: $\la
x-P_C(x), z-P_C(x)\ra\leq0$.
\end{proposition}
\proof See Theorem 3.14 of \cite{librobauch}. \qed

\noindent Next we deal with the so-called Fej\'er convergence and
its properties.
\begin{definition}\label{def-cuasi-fejer}
Let $S$ be a nonempty subset of $\RR^n$. A sequence
$(x^k)_{k\in\NN}$ is said to be Fej\'er convergent to $S$, if and
only if for all $x \in S$, there exists $\bar{k}>0$ such that $\|
x^{k+1}-x\| \leq \| x^{k}-x\|$ for all $k\ge\bar{k}$.
\end{definition}
This definition was introduced in \cite{browder} and has been
further elaborated in \cite{IST}. An useful result on Fej\'er
sequences is the following.
\begin{theorem}\label{cuasi-Fejer}
If $(x^k)_{k\in\NN}$ is Fej\'er convergent to $S$ then,
\item[i)] The sequence $(x^k)_{k\in\NN}$ is bounded,
\item [ii)] if a cluster point
of  the sequence $(x^k)_{k\in\NN}$ belongs to $S$, then the sequence
$(x^k)_{k\in\NN}$ converges to a point in $S$.
\end{theorem}
\proof  See Theorem $2.16$ of \cite{baus-borw}.\qed

\section{On K-convexity}\label{sec-3}

Convexity is a very helpful concept in optimization. Convex
functions satisfy nice properties such as existence of directional
derivative and subgradients, which are essential for optimality
conditions and iterative schemes for nonsmooth optimization
problems. In this section, we study the fulfillment of these
properties in the variable order case. First, we remind that
$F:\mathbb{R}^n\to\mathbb{R}^m$ is $K$-convex, respect to
$K:\RR^m\rightrightarrows \RR^m$ a point-to-cone application, if
\begin{equation}\label{convexidade1}
F(\alpha
x+(1-\alpha)\hat{x})\preceq_{K(F(\alpha x+(1-\alpha)\hat{x}))}\alpha
F(x)+(1-\alpha)F(\hat{x}),
\end{equation}
for any $x, \hat{x}\in\mathbb{R}^n$ and $\alpha\in[0,1]$ or equivalently \eqref{convexidade}.

\begin{remark} We want to point out that  this definition of convexity is independent of the concept introduced in \cite{yunier-gema}. There,  the condition is
\begin{equation}\label{eq.8}F(\alpha
x+(1-\alpha)\hat{x})\preceq_{K(\alpha x+(1-\alpha)\hat{x})}\alpha
F(x)+(1-\alpha)F(\hat{x}),\end{equation} for any $x, \hat{x}\in\mathbb{R}^n$ and $\alpha\in[0,1]$.
So, the order
is given by a  point-to-cone application $K$, whose domain is $\RR^n$ and not $\RR^m$ as in  \eqref{convexidade1}.
\end{remark}

Next examples show that there exist functions convex with respect to only one of two definitions presented in \eqref{convexidade1} and \cite{yunier-gema} (see \eqref{eq.8}).
\begin{example}
Let $F: \RR^2\to \RR^2$, $F(x_1,x_2)=(x^2_1+x_2^2+1, x_1)$, and $K:\RR^2\rightrightarrows\RR^2$,
$$K(x_1,x_2)=\left\{\begin{array}{ll}\RR^2_+,&\text{if} \ x_1\geq \frac{1}{2},
\\\\
\left\{r(\cos\theta\,,\,\sin\theta)\colon r\geq 0, \theta\in[\frac{3\pi}{4}-\frac{3\pi}{2}x_1,\frac{5\pi}{4}-\frac{3\pi}{2}x_1 ]\right\},&\text{if} \
x_1\in(0\,,\,\frac{1}{2}), \\\\
\left\{(z_1\,,\,z_2)\colon z_1\leq |z_2|\right\},  &\text{if} \
x_1\leq 0.\end{array}\right.$$ Note that $K(F(x))=\RR^2_+$ for all
$x\in \RR^n$. Since both components of $F$ are convex in the
classical sense, condition \eqref{convexidade1} holds and $F$ is
$K$-convex. However,
$$F(0,0)-\frac{F(x_1,x_1)+F(-x_1,-x_1)}{2}=(-2x_1^2\,,\,0)\notin -K(0,0), \text{ for all }x_1\neq0.$$
This means that the function is non-convex in the sense defined in \cite{yunier-gema} (see \eqref{eq.8}).\qed
\end{example}

\begin{example}
Let $F:[0,1]\times[0,1]\to\RR^2$, $F(x_1,x_2)=(x^2_1+x_2^2-5, x_2)$,  and $K:\RR^2\rightrightarrows\RR^2$,
$$K(x_1,x_2)=\left\{\begin{array}{ll}
\RR^2_+,&\text{if} \ x_1\geq -1,\\ \\
\left\{r(\cos\theta\,,\,\sin\theta)\colon r\geq 0,
\theta\in[-\pi-\pi x_1,-\frac{\pi}{2}-\pi x_1 ]\right\},& \text{if}
\ x_1\in (-2\,,\,-1), \\ \\ -\RR^2_+,  & \text{if} \ x_1\leq -2.
\end{array}\right.$$
Actually, for all $x$ belonging to the domain of $F$, i.e., the set
$[0,1]\times[0,1]$, $K(x)=\RR^2_+$ and so, $F$ is  convex with
respect to the order defined in \cite{yunier-gema} (see
\eqref{eq.8}). That is,  for all $x,\hat{x}\in [0,1]\times[0,1]$
$$F(\alpha
x+(1-\alpha)\hat{x})\preceq_{K(\alpha x+(1-\alpha)\hat{x})}\alpha
F(x)+(1-\alpha)F(\hat{x}).$$
On the other hand, the image of $F$ lies in $[-5,-3]\times[0,1]$, This means that $K(F(x))=-\RR^2_+$. Since
the vector
$F(\alpha x+(1-\alpha)\hat{x})-\alpha F(x)-(1-\alpha)F(\hat{x})$
is not identically $0$ and, as already remarked,  it belongs to $\RR^2_+$ for all $x,\hat{x}\in [0,1]\times[0,1]$, \eqref{convexidade1} is not fulfilled.\qed
\end{example}

Now we begin with the analysis of the $K$-convexity defined in \eqref{convexidade1}. First the epigraph of $K$-convex functions will be studied.
In the variable order case the epigraph of $F$ is defined  as
$$epi(F):=\{(x,y)\in \RR^n\times \RR^m:\; F(x)\in y-K(F(x))\}.$$
In non-variable orders, \emph{i.e.}, when $K$ is a constant
application, the convexity of $epi(F)$ is equivalent to the
convexity of $F$; see \cite{luc}. However, as it is shown in the
next proposition, in the variable order setting this important
characterization does not hold.

\begin{proposition} Suppose that  $F$ is a $K$-convex function. Then, $epi(F)$ is convex if and only if $K(F(x))\equiv K$, for all $x\in \RR^n$.
\end{proposition}
\proof Suppose that  for some $x,\hat{x}\in \RR^n$ such that
$F(x)\neq F(\hat{x})$, there exists $z\in K(F(x))\setminus
K(F(\hat{x}))$. Take the points  $\left(x, F(x)+2\alpha z\right)$
and $\left(2\hat{x}-x,F(2\hat{x}-x)\right)$, with $\alpha>0$. They
belong to $epi(F)$.

\noindent Consider the following convex combination:  $$\frac{(x,
F(x)+2\alpha
z)}{2}+\frac{(2\hat{x}-x,F(2\hat{x}-x))}{2}=\left(\hat{x},\frac{F(x)+F(2\hat{x}-x)}{2}+
\alpha z\right).$$ This point belongs to $epi(F)$ if and only if
$$F(\hat{x})= \frac{F(x)+F(2\hat{x}-x)}{2}+\alpha z- k(\alpha),$$ where  $k(\alpha)\in K(F(\hat{x}))$.
By the $K$-convexity of $F$,
$$F(\hat{x})=\frac{F(x)+F(2\hat{x}-x)}{2}-k_1,$$ where $k_1\in
K(F(\hat{x}))$. So, \begin{equation}\label{hiperplano1}\alpha z+
k_1= k(\alpha).\end{equation} Since $K(F(\hat{x}))$ is closed and
convex, and $z\notin K(F(\hat{x}))$, $\{z\}$ and $K(F(\hat{x}))$ may
be strictly separated in $\RR^m$ by a hyperplane, \emph{i.e.}, there
exists some $p\in \RR^m\setminus \{0\}$ such that
\begin{equation}\label{hiperplano}p^Tk\geq 0> p^Tz,\end{equation}
for all
 $k\in K(F(\hat{x}))$. Therefore,
 after multiplying \eqref{hiperplano1} by $p^T$   and
using \eqref{hiperplano}  with $$k=k(\alpha)\in K(F(\hat{x})),$$ we obtain that
\begin{equation*}
 \alpha p^T z+ p^T k_1= p^Tk(\alpha)\geq 0.
\end{equation*}
Taking limits as  $\alpha$ goes to $\infty$, the contradiction is
established, because  $$0\leq \alpha p^T z+p^Tk_1\to -\infty .$$
Hence,
 $K(F(x) )\equiv K$ for all $x\in \RR^n. $
\qed

In the following we present some analytical properties of $K$-convex
functions. For the  non-differentiable model, we generalize the
classical assumptions given in the case of constant cones; see
\cite{jan,luc}. Let us first present the definition of Daniell cone,
for more details; see \cite{Peressini}.

Let $\mathcal{K}$ be a closed and convex cone. Given the partial
order structure induced by a cone $\mathcal{K}$, the concept of
infimum of a sequence can be defined. Indeed, for a sequence
$(x^k)_{k\in \NN}$ and a cone $\mathcal{K}$, the point $\hat{x}$ is
$\dsty\inf_{k\in \NN} \{x^k\}$ if and only if $(x^k-\hat{x})_{k\in
\NN}\subset \mathcal{K}$, and there is not $\bar{x}$ such that
$\hat{x}-\bar{x}\in \mathcal{K}$ and $(x^k-\bar{x})_{k\in\NN}\subset
\mathcal{K}$.
\begin{definition} We say that a convex cone $\mathcal{K}$ is Daniell cone
iff, for all sequence $(x^k)_{k\in \NN}\subset\RR^n$ satisfying
$(x^k-x^{k+1})_{k\in \NN}\subset \mathcal{K}$ and for some $x\in
\mathbb{R}^n$,  $(x^k-x)_{k\in \NN}\subset \mathcal{K}$, then
$\displaystyle\lim_{k\rightarrow\infty}x^k=\inf_{k\in
\NN}\{x^k\}$.\end{definition}

It is well known that  every pointed, closed and convex cone in a
finite dimensional space is a Daniell cone; see, for instance,
\cite{theluc}.

\begin{lemma}\label{lemmita} Suppose that there exists $\mathcal{K}$ a Daniell cone such that $K(F(x)) \subseteq \mathcal{K}$
for all $x$ in a neighborhood of $\hat{x}$.
If $F$ is a $K$-convex function,  then $F$ is locally Lipschitz around $\hat{x}$.
\end{lemma}
\proof If $F$ is $K$-convex, then $F$ is $\mathcal{K}$-convex in the
non-variable sense. By Theorem $3.1$ of \cite{luc-subdif}, $F$ is
locally Lipschitz.\qed

\begin{proposition}
Suppose that for each $\bar{x}$ there exists $\varepsilon>0$ such
that $\cup_{x\in B[\bar{x},\varepsilon]} K(F(x))\subseteq
\mathcal{K}$, where $\mathcal{K}$ is a Daniell cone. Then, the
directional derivative of $F$ at $\bar{x}$ exists along
$d=x-\bar{x}$, that is,
$$F^\prime(\bar{x};x-\bar{x})=\lim_{t\to 0^+} \frac{F(\bar{x}+td) -F(\bar{x})}{t}.$$
\end{proposition}
\proof By the convexity of $F$,
$$F(\bar{x}+t_1d)-\frac{t_1}{t_2}F(\bar{x}+t_2d)-\left(\frac{t_2-t_1}{t_2}\right)F(\bar{x})\in
-K(F(\bar{x}+t_1d)),$$ for all $0<t_1<t_2<\varepsilon.$ Dividing by
$t_1$, we have $$\frac{F(\bar{x}+t_1d)-F(\bar{x})}{t_1}-
\frac{F(\bar{x}+t_2d)-F(\bar{x})}{t_2} \in
-K(F(\bar{x}+t_1d))\subseteq -\mathcal{K}.$$
\\Hence, $\displaystyle \frac{F(\bar{x}+t_1d)-F(\bar{x})}{t_1}$ is a non-increasing function.
Similarly, as
$$F(\bar{x})-\frac{t_1}{t_1+1}F(\bar{x}-d)-\frac{1}{t_1+1}F(\bar{x}+t_1d)\in -K(F(\bar{x})),$$
it holds that   $$\frac{F(\bar{x}+t_1d)-F(\bar{x})}{t_1}-
F(\bar{x}-d)-F(\bar{x})\in K(F(\bar{x}))\subseteq \mathcal{K}.$$
Since $\mathcal{K}$ is a Daniell cone,
$\displaystyle\frac{F(\bar{x}+t_1d)-F(\bar{x})}{t_1}$ has a limit as
$t_1$ goes to $0$. Hence, the directional derivative exists. \qed

\noindent Let us present the definition of subgradient.

\begin{definition} We say that
 $\epsilon_{\bar{x}}\in \RR^{m\times n}$ is a subgradient of $F$ at $\bar{x}$ if for all $x\in \RR^n$, $$ F(x)-F(\bar{x})\in \epsilon_{\bar{x}} (x-\bar{x})+K(F(\bar{x})).$$
The set of all subgradients of $F$ at $\bar{x}$ is denoted as $\partial F(\bar{x})$.
\end{definition}

\begin{proposition}
 If for all $x\in \RR^n$,  $\partial F(x)\neq \emptyset$, then $F$ is $K$-convex.
\end{proposition}

\proof Since $\partial F(x)\neq \emptyset$, for all $x\in \RR^n$,
taking any $\bar{x}, \hat{x}\in\RR^n$ there exists $\epsilon_{\alpha
\bar{x}+(1-\alpha)\hat{x}}$ belonging to $\partial F(\alpha
\bar{x}+(1-\alpha)\hat{x})$ and $k_1,k_2\in K(F(\alpha
\bar{x}+(1-\alpha)\hat{x}))$, such that
$$F(\hat{x})-F(\alpha \bar{x}+(1-\alpha)\hat{x})=\alpha\epsilon_{\alpha
\bar{x}+(1-\alpha)\hat{x}}(\hat{x}-\bar{x}) +k_1, $$and
$$F(\bar{x})-F(\alpha \bar{x}+(1-\alpha)\hat{x})=(\alpha-1)\epsilon_{\alpha \bar{x}+(1-\alpha)\hat{x}}(\hat{x}-\bar{x}) +k_2.$$
Multiplying the previous equalities by $(1-\alpha)$ and $\alpha$ respectively, their addition leads to
$$\alpha F(\bar{x})+(1-\alpha)F(\hat{x})-F(\alpha \bar{x}+(1-\alpha)\hat{x})= \alpha k_2+ (1-\alpha) k_1.$$
Since $K(F(\alpha \bar{x}+(1-\alpha) \hat{x}))$ is convex, the
result follows. \qed

\begin{proposition}\label{prop11}If $K$ is a closed application, then $\partial F$ is closed.\end{proposition}
\proof Assume that $(x^k)_{k\in\NN}$ and $(A^k)_{k\in\NN}$ are
sequences such that $A^k\in
\partial F(x^k)$ for all $k$, $\lim_{k\to\infty} x^k= \bar{x}$ and
$\lim_{k\to\infty} A^k= A$. For every $x$, one has
$$F(x) - F(x^k)- A^k(x - x^k) \in K(F(x^k)).$$
Taking $k$ going to $\infty$,  as $\lim_{k\to\infty}
F(x^k)=F(\bar{x})$ and $K$ is a closed mapping, we get that
$$F(x) - F(\bar{x})- A(x - \bar{x})\in K(F(\bar{x})).$$
Hence, $A\in \partial F(\bar{x})$, establishing that $\partial
F(\bar{x})$ is closed.\qed

\begin{proposition}\label{propgrad}
 Let $F$ be a $K$-convex function.
 If $gr(K)$ is closed, then for all $\bar{x}\in \RR^n$, where $F$ is differentiable, $\nabla F(\bar{x})=\partial F(\bar{x})$.
\end{proposition}

\proof
First we  show that $\nabla F(\bar{x})$ belongs to  $\partial F(\bar{x})$.
  Since $F$ is a differentiable function, fixed $\bar{x}$, we get
 $$F(\alpha x+(1-\alpha) \bar{x})=F(\bar{x})+\alpha\nabla F(\bar{x})(x-\bar{x}) +o\left(\alpha\right).$$
By $K$-convexity,
 $$F(\bar{x})+\alpha\nabla F(\bar{x})(x-\bar{x}) +o\left(\alpha\right)\in \alpha F(x)+(1-\alpha)F(\bar{x})-K(F(\alpha x +(1-\alpha) \bar{x})).$$
So,
$$  \alpha\left( F(x)- F(\bar{x})-\nabla F(\bar{x})(x-\bar{x}) + \frac{o(\alpha)}{\alpha}\right) \in K(F(\alpha x +(1-\alpha) \bar{x})).$$
Since $K$ is a cone, it follows that
$$   F(x)- F(\bar{x})-\nabla F(\bar{x})(x-\bar{x}) +\frac{o(\alpha)}{\alpha}\in K(F(\alpha x +(1-\alpha) \bar{x})).$$
By taking   limits as $\alpha$ goes to $0$ and recalling that $F$ is a continuous function and $K$ is a closed application, by Lemma \ref{lemmita}  it holds that

$$   F(x)- F(\bar{x})-\nabla F(\bar{x})(x-\bar{x}) \in K(F( \bar{x}) ),$$
and hence,  $\nabla F(\bar{x})\in \partial F(\bar{x})$.

\noindent Suppose that $\varepsilon_{\bar{x}}\in \partial
F(\bar{x})$. Fixed $d\in \RR^n$, we get that, for all $\alpha> 0$,
$$F(\bar{x} +\alpha d)-F(\bar{x})=\alpha \nabla F(\bar{x})d+o(\alpha)\in \alpha\, \varepsilon_{\bar{x}} d+k(\alpha),$$
where $k(\alpha)\in K(F(\bar{x}))$. Dividing by $\alpha>0$, and
taking limits as $\alpha$ approaches  $0$, it follows that
$$[\nabla F(\bar{x})-\varepsilon_{\bar{x}}]d\in K(F(\bar{x})),$$ recall that $K(F(\bar{x}))$ is a closed set. Repeating the same analysis for $-d$, we obtain that
 $$-[\nabla F(\bar{x})-\varepsilon_{\bar{x}}]d\in K(F(\bar{x})).$$
Taking into account that $K(F(\bar{x}))$ is a pointed cone, $[\nabla
F(\bar{x})-\varepsilon_{\bar{x}}]d=0$. As the previous equality is
valid for all $d\in \RR^n$,
$$\nabla F(\bar{x})=\varepsilon_{\bar{x}},$$ establishing the desired equality.
\qed

\begin{theorem}\label{t2}Suppose that there exists $\mathcal{K}$ a Daniell cone such that $K(F(x)) \subseteq \mathcal{K}$ for all $x$ in a neighborhood of $\hat{x}$.
If $F$ is $K$-convex  and $K$ is a  closed application, then $\partial F(\hat{x})\neq\emptyset$.
\end{theorem}
\proof By Lemma \ref{lemmita}, $F$ is a locally Lipschitz continuous function. By Rademacher's Theorem, for all $\hat{x}$,
 $F$ is differentiable almost everywhere on some neighborhood of
$\hat{x}$. Moreover, due to the boundedness of $\nabla F$ whenever exists,  there exists a sequence $x^k$ convergent to $\hat{x}$ such that
$A=\lim_{k\rightarrow \infty} \nabla F(x^k)$.
By Proposition \ref{propgrad}, it holds that  $\nabla F(x^k)= \partial F(x^k)$. By Proposition \ref{prop11}, $A\in \partial F(\hat{x})$,
hence $\partial F(\hat{x})\neq \emptyset$. \qed
\begin{remark}\label{nota1} Given  $\hat{x}$ and $V$ a bounded neighborhood of $\hat{x}$, under the assumptions of the previous Theorem, the set
$\partial F(x)$ is uniformly bounded in $ V$. Indeed as $F$ is
$K$-convex, locally around $\hat{x}$, $F$ will be also
$\mathcal{K}$-convex. Now, since the domain of $F$ is a finite
dimensional space, the fact follows directly by \cite[Theorem
4.12(ii)]{luc-subdif}.
 \end{remark}

\section{The Algorithms}\label{sec-4}
In this section we consider two variants of subgradient method for
solving problem  \eqref{prob-1}. The algorithms generate a sequence
of projections onto special sets. From now on, we assume that the
following assumptions hold.

\medskip

\noindent {\bf{Assumptions}}
\vspace{-5pt}
\begin{itemize}\item [(A1)] The subgradients of $F$ are locally bounded.
 \item[(A2)] $F$ is $K$-convex.\item[(A3)] $K:\RR^m\rightrightarrows\RR^m$ is a closed application.
 \item[(A4)] For all $x^*\in S_*$ and $x\in C,$
 \begin{equation}\label{cond-cone}
K(F(x^*))\subseteq K(F(x)).\end{equation}
\end{itemize}
We emphasize that Assumption (A1) is a typical hypothesis for
proving the convergence of the subgradient-scalar methods in
infinite dimension setting; see \cite{sc-1, AIS, poljak,
yunier-iusem-1}. As stated in \cite{luc-subdif}, for the scalar and
vector framework, this assumption holds trivially in
finite-dimensional spaces. Recently, (A1) was proved in
\cite{yun-sub}, when $K$ is a constant application. A sufficient
condition can be found in Remark \ref{nota1}.

\noindent The existence of subgradient is guaranteed in Theorem \ref{t2}.

Assumption (A4) implies that there exists a cone $\mathbb{K}$ such that  $K(F(x^*))\equiv  \mathbb{K}$  for all $x^*\in S_*$.
In  this  case problem  \eqref{prob-1}
is equivalent to the non-variable inequalities problem
$$\text{find $x\in C$ such that }F(x)\preceq_{\mathbb{K}} 0.$$
However, as $\mathbb{K}$ is not known, this equivalence is not useful from a practical viewpoint. Next example shows a function and an order structure
fulfilling
\eqref{cond-cone}.

\begin{remark} Given problem \eqref{prob-1} with $C=\RR$, $F:\RR\to \RR^2$, $F(x)=(x^2,x)$,
$K\colon\RR^2\rightrightarrows\RR^2$, $K(y)=\{r(\cos\theta\,,\,\sin\theta)\colon r\geq 0, \theta\in[0,\theta(y)]\}$, where
$$\theta(y)=\left\{\begin{array}{ll}
\displaystyle \frac{\pi}{2},&\text{ if }y_1=0,\\ \\  \displaystyle \frac{3\pi}{4}-\frac{\arctan(y_2^2/y_1^2)}{2}, &\text{ otherwise. }
\end{array}\right.$$
Evidently $$\RR_+\times\{0\}\subset K(y)\subset\RR_+\times \RR.$$
Moreover, $F(x)\in -K(F(x))$ if and only if $x=0$. Therefore,
$S_*=\{0\}$ and due to
$$\theta(y)\geq \frac{\pi}{2}=\theta(0,0),$$ Assumption (A4) holds.

\noindent Since $F_1(x)=x^2$ is convex and $F_2(x)=x$ is a linear function,
$$F(\alpha x+(1-\alpha)\hat{x})-\alpha F(x)-(1-\alpha) F(\hat{x})\in-\RR_+\times\{0\}\subseteq -K(F(\hat{x}))$$
 for all $x, \hat{x}\in \RR$. Hence, $F$ is $K$-convex. Moreover, the continuity of $\theta$ implies that $K$ is a closed application.
\end{remark}
\noindent Now we will present the conceptual algorithm.

\vspace{10pt}

\noindent {\bf Conceptual Algorithm}

\medskip
\noindent{\bf Initialization step.} Take $x^0\in C$, and set $k=0$.

\medskip

\noindent{\bf Iterative step.} Given $x^k$, $U^k\in\partial F(x^k)$. Compute
\begin{equation}\label{paso3_*}
\displaystyle x^{k+1}:=\mathcal{F}(x^k, U^k).
\end{equation}

\noindent If $x^{k+1}=x^k$ then stop.

\medskip

\noindent We consider two variants of the conceptual algorithm.  As
they are based on the algorithms proposed in \cite{rob,yun-luis},
the extensions are called Algorithms $R$ and $S$ respectively. The
main difference is given by the definition of the procedure
$\mathcal{F}$ in \eqref{paso3_*}, which is defined as follows

\begin{eqnarray}\label{Alg-R}
\mathcal{F}_R(x^k, U^k) &:=& P_{C\cap H(x^k,U^k)}(x^k);\\
\label{Alg-S} \mathcal{F}_S(x^k, U^k) &:=& P_{C\cap W(x^k)\cap
H(x^k,U^k)}(x^0);\end{eqnarray}
 where
\begin{equation*}\label{HHH}H(x,U):=\{z\in \RR^n: F(x)+U(z-x)\in -K(F(x))\}\end{equation*} and
\begin{equation*}\label{Wk}
W(x):=\left\{z\in \RR^n\,:\,\la z-x,x^0-x\ra\leq 0\right\}.
\end{equation*}
Before we start with the formal analysis of the convergence
properties of the algorithm, we make a comment on the complexity of
the projection steps, defined in \eqref{Alg-R} and \eqref{Alg-S}.
First, we want to point out that $W(x)$ is a halfspace and $H(x,U)$
is convex by the convexity of $-K(F(x))$ for any $x\in C$.
Furthermore, if the dual cone of $K(F(x))$,
$$K^*(F(x)):=\left\{z\in \RR^m\colon \langle z\,,\, y\rangle\geq 0, \; \forall y\in K(F(x)) \right\},$$ has finitely many generators, that is,
exist $G=\{u_1,u_2,\ldots,u_r\}\subset K^*(F(x))$, such that
$$K^*(F(x))=\left\{z\in \RR^m\colon z=\sum_{i=1}^r \lambda_iu_i, \lambda_i\geq 0, \; i=1, \ldots, r \right\},$$
then $ H(x,U)$ is the intersection of $r$ halfspaces.

\begin{remark}
Note that, if $C$ is described by nonlinear constrains, the addition
of linear constraints may lead to a smaller set, onto which it may
be easier to project;  see, for instance, \cite{yunier-iusem-1}. So,
if  $ K^*(F(x^k))$ has finitely many generators,  the sets
$H(x^k,U^k)$ and $H(x^k,U^k)\cap W(x^k)$ are the intersection of
finitely many halfspaces, as was noted above. Thus, the projections
defined in \eqref{Alg-R} and \eqref{Alg-S} do not entail any
significant additional computational cost over the computation of
the projection onto $C$ itself.
\end{remark}

\section{Convergence Analysis}\label{sec-5}
In this part we prove the convergence of the algorithms. The section
will contain three subsections. First we study the properties of the
solution set $S_*$ and present some general properties of the
conceptual algorithm. The convergence analysis of the proposed
variants, Algorithms $R$ and $S$, will be presented separately in
the last two subsections.

\subsection{Properties of the Solution Set}
\begin{proposition}\label{cerradita}The set $S_*$ is  closed and convex.
\end{proposition}
\proof Take $x,x^{*}\in S_*$. Then, it holds that

$$F(\alpha x+(1-\alpha)x^*)\in \alpha F(x)+(1-\alpha)F(x^*)- K(F(\alpha x+(1-\alpha)x^*)),$$ for all $\alpha\in [0,1]$.
Since $F(x) \preceq_{K(F(x))} 0$ and   $F(x^*) \preceq_{K(F(x^*))}
0$, it follows from (A4) that $$K(F(x))=K(F(x^*))\subseteq
K(F(\alpha x+(1-\alpha)x^*)) .$$ Hence, $$F(\alpha
x+(1-\alpha)x^*)\in - K(F(\alpha x+(1-\alpha)x^*)),$$ and therefore
$\alpha x+(1-\alpha)x^*\in S_*$.

For the closeness, consider any sequence $(x^k)_{k\in\NN}\subset
S_*$ convergent to  $x^*$. Since $F$ is a continuous function; see
Lemma \ref{lemmita}, $\lim_{k\rightarrow\infty} F(x^k)=F(x^*)$ and
taking into account that $F(x^k)\in -K(F(x^k))$ and the closedness
of $K$ leads to $F(x^*)\in -K(F(x^*))$. So, $x^*\in S_*$. \qed

\noindent We assume that $S_*$ is a nonempty set.

\begin{lemma}\label{sol-en-H}For all $x \in C\setminus S_*$  and $U\in\partial F(x)$, it holds that  $S_*\subseteq H(x, U)$.
\end{lemma}
\proof Take $x^*\in S_*$. Then, $F(x^*)\in -K(F(z))$ and by the subgradient inequality,
 $$F(x)+U(x^*-x)-F(x^*)\in -K(F(x)),$$ for all $x\in C$ and all $U\in \partial F(x)$.
 Hence, using the above inclusion and \eqref{cond-cone}, we get that
  $$F(x)+U(x^*-x)\in -K(F(x)) - K(F(x^*))\subseteq  -K(F(x)),$$ for all $x\notin S_*$. So, $x^*\in H(x, U)$.  \qed
 \begin{lemma}\label{lemma3}  If $x\in H(x, U)\cap C$ for some $U\in\partial F(x)$, then $x\in S_*$.
\end{lemma}
\proof Suppose that $x\in H(x, U)\cap C$ for some $U\in\partial
F(x)$, then $x\in C$  and $$F(x)\in -K(F(x)),$$ \emph{i.e.},  $x\in
S_*$.  \qed

The above lemma will be useful to show that the stop criterion of
the variants of the conceptual algorithm are well defined.
\subsection{Convergence of Algorithm R}
In this subsection all results are referent to Algorithm $R$, i.e.,
with the iterative step as
$$x^{k+1}=\mathcal{F}_R(x^k, U^k) = P_{C\cap H(x^k,U^k)}(x^k),$$
where \begin{equation*}H(x^k,U^k)=\{z\in \RR^n: F(x^k)+U^k(z-x^k)\in -K(F(x^k))\}\end{equation*} and $U^k\in\partial
F(x^k)$.

\noindent The following proposition gives the validity of the stop
criterion on Algorithm $R$.
\begin{proposition}
If Algorithm $R$ stops at iteration $k$, then $x^k\in S_*$.
\end{proposition}
\proof If Algorithm $R$ stops, then $x^{k+1}=x^k$. It follows from
\eqref{Alg-R} that $x^k\in H(x^k,U^k)\cap C$.  So, by Lemma
\ref{lemma3}, $x^k\in S_*$. \qed
\begin{proposition}\label{AlgR-fejer}The sequence generated by Algorithm $R$ is F\'ejer convergent
 to $S_*$. Moreover, it is bounded and \begin{equation*}\label{oppRob}
 \lim_{k\rightarrow\infty}\|x^{k+1}-x^k\|= 0.\end{equation*}
\end{proposition}
\proof Take $x^*\in S_*$. By Lemma \ref{sol-en-H},  $x^*\in  H(x^k,
U^k)$, for all $k\in \mathbb{N}$. Then
$$\|x^{k+1}-x^*\|^2-\|x^{k}-x^*\|^2
+\|x^{k+1}-x^k\|^2=2\langle x^*-x^{k+1}, x^k-x^{k+1}\rangle\leq
0,$$ using Proposition \ref{popiedades_projeccion} and \eqref{Alg-R}
in the last inequality. So,
\begin{equation}\label{des-fejer-algR}\|x^{k+1}-x^*\|^2\leq
\|x^{k}-x^*\|^2-\|x^{k+1}-x^k\|^2.\end{equation} The above
inequality implies that $(x^{k})_{k\in\NN}$ is Fej\'er convergent to
$S_*$ and hence $(x^k)_{k\in\NN}$ is bounded. We get
$$0\leq \|x^{k+1}-x^*\|^2\leq \|x^{k}-x^*\|^2.$$
So, $(\|x^{k}-x^*\|^2)_{k\in\NN}$ is a convergent sequence.
Therefore, using \eqref{des-fejer-algR}, we obtain that

\begin{equation*}
 \lim_{k\rightarrow\infty}\|x^{k+1}-x^k\|= 0.
\end{equation*}
\qed
\begin{theorem}\label{conv-Alg-R}
The sequence generated by Algorithm $R$ converges to some point in $S_*$.
\end{theorem}
\proof By Proposition \ref{AlgR-fejer}, $(x^k)_{k\in\NN}$ is
bounded. So, using (A1), $(U^k)_{k\in\NN}$ is bounded, \emph{i.e.},
there exists $L\ge 0$ such that
\begin{equation}\label{op-p}
 \|U^k\|\leq L,
\end{equation} for all $k$.

Fix $k\in \mathbb{N}$. Since $K(F(x^k))$ is a closed convex cone, it is clear that   $y\in \RR^m$ can be uniquely written  as
$$y=y_++y_-,$$ with $y_+\in K^*(F(x^k))$,  $y_-\in -K(F(x^k))$ and $\langle y_+, y_-\rangle=0$.
For $y=F(x^k)$, consider $F(x^k)_+$ and $F(x^k)_-$.
Now
\begin{eqnarray*}\|F(x^k)_+\|^2&=&\left\langle F(x^k)_+,F(x^k)_+ +F(x^k)_-\right\rangle=\left\langle F(x^k)_+,F(x^k)\right\rangle \\
&=& \left\langle F(x^k)_+,F(x^k) +U^k(x^{k+1}-x^k)\right\rangle
-\left\langle F(x^k)_+,U^k(x^{k+1}-x^k)
\right\rangle.\end{eqnarray*}

\noindent But $F(x^k)_+\in K^*(F(x^k))$, so $\left\langle
F(x^k)_+,F(x^k) +U^k(x^{k+1}-x^k)\right\rangle \leq 0 $ and,
therefore
\begin{equation*}\|F(x^k)_+\|^2\leq -\left\langle F(x^k)_+,U^k(x^{k+1}-x^k) \right\rangle.\end{equation*}
Applying the Cauchy Schwartz inequality and recalling \eqref{op-p}, it follows that
\begin{equation*}
 \|F(x^k)_+\|^2\leq L\|F(x^k)_+\|\|x^{k+1}-x^k\|.
\end{equation*}
Since $x^{k}\notin S_*$, $F(x^k)_+\neq 0$. So, dividing by $
\|F(x^k)_+\|$, we obtain
\begin{equation*}
 \|F(x^k)_+\|\leq L\|x^{k+1}-x^k\|.
\end{equation*}
Recalling Proposition \ref{AlgR-fejer},  it follows that
\begin{equation}\label{oppp-p}
\lim_{k\rightarrow\infty} \|F(x^k)_+\|= 0.
\end{equation}
Now consider a convergent subsequence $(x^{\ell_k})_{k\in\NN}$ of
$(x^k)_{k\in\NN}$. Denote $x^*$ as its limit. It follows from
\eqref{oppp-p} that $F(x^*)_+=0$.  Henceforth, $F(x^*)=F(x^*)_-$.
Moreover as
$$\lim_{k\to\infty}F(x^{\ell_k})_-=\lim_{k\to\infty}F(x^{\ell_k})-\lim_{k\to\infty}F(x^{\ell_k})_+,$$
we get that
$$\lim_{k\to\infty}F(x^{\ell_k})_-=F(x^{*}).$$
Since $F(x^{\ell_k})_-\in -K(F(x^{\ell_k}))$ and  (A3) is fulfilled,
$$F(x^{*})\in -K(F(x^{*})),$$
\emph{i.e.}, $x^{*}\in S_*$. Therefore, the accumulation points of
$(x^k)_{k\in\NN}$ belong to $S_*$. Finally, by the F\'ejer
convergence, the sequence converge to a point in $S_*$.\qed

\subsection{Convergence of Algorithm S}
In this subsection all results are referent to Algorithm $S$, i.e.,
with the iterative step as
$$x^{k+1}=\mathcal{F}_S(x^k, U^k) = P_{C\cap W(x^k)\cap H(x^k,U^k)}(x^0),$$
where \begin{equation}\label{*}H(x^k,U^k)=\left\{z\in \RR^n:
F(x^k)+U^k(z-x^k)\in -K(F(x^k))\right\}\end{equation} with
$U^k\in\partial F(x^k)$ and \begin{equation}\label{**}
W(x^k)=\left\{z\in \RR^n\,:\,\la z-x^k,x^0-x^k\ra\leq 0\right\}.
\end{equation}
The following proposition gives the validity of the stop criterion
on Algorithm $S$.
\begin{proposition}
If Algorithm $S$ stops at iteration $k$, then $x^k\in S_*$.
\end{proposition}
\proof If Algorithm $S$ stops at iteration $k$, then $x^{k+1}=x^k$.
It follows from \eqref{Alg-S} that $x^k\in W(x^k)\cap H(x^k,U^k)\cap
C\subseteq H(x^k,U^k)\cap C$. So, by Lemma \ref{lemma3}, $x^k\in
S_*$. \qed

Observe that, in virtue of their definitions, given in \eqref{*} and
\eqref{**}, $W(x^k)$ and $H(x^k,U^k)$ for some $U^k\in\partial
F(x^k)$ are convex and closed sets, for each $k\in \NN$. Therefore
$C\cap H(x^k,U^k)\cap W(x^k)$ is a convex and closed set, for each
$k\in \NN$. So, if $C\cap H(x^k,U^k)\cap W(x^k)$ is nonempty then,
the next iterate, $x^{k+1}$, is well-defined. Next lemma guarantees
this fact.

\begin{lemma}\label{lemma:3} For all $k\in \NN$, it holds that $S_*\subseteq C\cap H(x^k,U^k)\cap W(x^k)$.
\end{lemma}
\proof We proceed by induction. By definition, $S_*\subseteq C$. By
Lemma \ref{sol-en-H}, $S_*\subseteq C\cap H(x^k,U^k)$, for all $k$.
For $k=0$, since $W(x^0)=\RR^m$, $S_*\subseteq C\cap H(x^0,U^0)\cap
W(x^0)$. Assume that $S_*\subseteq C\cap H(x^{\ell},U^{\ell}) \cap
W(x^{\ell})$, for all $0\le\ell\leq k$. Henceforth,
$x^{k+1}=P_{C\cap H(x^k,U^k)\cap W(x^k)}(x^0)$ is well defined.
Then, by Lemma \ref{sol-en-H}, for all $x^*\in S_*$, we get that
\begin{equation*}\label{x*ink+1}\left\langle x^*-x^{k+1}\,,\, x^0-x^{k+1}\right\rangle=\left \langle x^*-P_{C\cap H(x^k,U^k)\cap W(x^k)}(x^0)\,,\,
x^0-P_{C\cap H(x^k,U^k)\cap
W(x^k)}(x^0)\right\rangle\leq0,\end{equation*} using the induction
hypothesis. The above inequality implies that $x^*\in W(x^{k+1})$
and hence, $S_*$ is a subset of $ C\cap H(x^{k+1},U^{k+1})\cap
W(x^{k+1})$. \qed

\begin{corollary}\label{l:well-definedness}  Algorithm $S$ is well-defined.
\end{corollary}
\proof By the previous lemma, $ S_*\subseteq C\cap H(x^k,U^k)\cap
W(x^k)$, for $k\in \NN$. Since  $S_*\neq\emptyset$, then,  given
$x^0$, the sequence $(x^k)_{k\in\NN}$ is computable. \qed

Before proving the convergence of the sequence, we will study its
boundedness. Next lemma  shows that  the sequence remains in a ball
determined by the initial point.

\begin{lemma}\label{l:limitacao} The sequence $(x^k)_{k\in\NN}$ is bounded. Furthermore,
\begin{equation*}\label{eq:bolas}
(x^k)_{k\in\NN}\subset
B\left[\frac{1}{2}(x^0+x^*),\frac{1}{2}\rho\right],
\end{equation*}
where $x^*=P_{S_*}(x^0)$ and $\rho={\rm dist}(x^0, S_*)$.
\end{lemma}
\proof Lemma \ref{lemma:3} says that $S_*\subseteq C\cap W(x^k)\cap
H(x^k,U^k)$ for   $k\in \NN$ and, by the definition of $x^{k+1}$ in
\eqref{paso3_*} and \eqref{Alg-S}, it is true that
\begin{equation}\label{eq:12}
\| x^{k+1}-x^0\| \leq\| z-x^0\|,
\end{equation}
for  $k\in \NN$ and all $z\in S_*$. Henceforth, taking in
\eqref{eq:12} $z=x^*$,
\begin{equation*}\label{eq:yunier}
 \| x^{k+1}-x^0\|\leq\|x^*-x^0\|=\rho,
\end{equation*}
for all $k$. Hence, $(x^k)_{k\in\NN}$ is bounded. Without loss of
generality, take $z^{k}=x^{k}-\frac{1}{2}(x^0+x^*)$ and
$z^{*}=x^{*}-\frac{1}{2}(x^0+x^*)$. It follows from the fact
$x^{*}\in W(x^{k+1})$ that
\begin{eqnarray*}
0&\geq& 2\la x^*-x^{k+1},x^0-x^{k+1}\ra \\&=&2\left\la
z^*+ \frac{1}{2}(x^0+x^*)-z^{k+1}-\frac{1}{2}(x^0+x^*),z^0+\frac{1}{2}(x^0+x^*)-z^{k+1}-\frac{1}{2}(x^0+x^*)\right\ra\\&=&2\left\la
z^*-z^{k+1},z^0-z^{k+1}\right\ra=\left\la
z^*-z^{k+1},-z^*-z^{k+1}\right\ra = \|z^{k+1}\|^2-\|z^*\|^2,
\end{eqnarray*}
using in the third equality that $z^*=-z^0$. So,
\begin{equation*}\label{eq:raio}
\left
\|x^{k+1}-\frac{x^0+x^*}{2}\right\|\leq\left\|x^*-\frac{x^0+x^*}{2}\right\|=\frac{\rho}{2},
\end{equation*} establishing the result.
\qed

\noindent  Now we will focus on the properties of the accumulation points.

\begin{lemma}\label{l:optimalidad} All accumulation points of $(x^k)_{k\in\NN}$ are elements of $S_*$.
\end{lemma}

\proof Since $x^{k+1}\in W(x^k)$,
\begin{equation*}
0\geq 2 \la
x^{k+1}-x^k,x^0-x^k\ra=\|x^{k+1}-x^k\|^2-\|x^{k+1}-x^0\|^2+\|x^k-x^0\|^2.
\end{equation*}
Equivalently,
$$0\leq\|x^{k+1}-x^k\|^2\leq\|x^{k+1}-x^0\|^2-\|x^k-x^0\|^2,$$
establishing that $(\|x^k-x^0\|)_{k\in\NN}$ is a monotone
nondecreasing sequence. It follows from Lemma \ref{l:limitacao} that
$(\|x^k-x^0\|)_{k\in\NN}$ is bounded and thus, it is a convergent
sequence. Therefore,
\begin{equation*}\label{xk+1-xk-va-cero}
\lim_{k\rightarrow\infty}\| x^{k+1}-x^k\|=0.
\end{equation*}
Let $\bar{x}$ be an accumulation point of $(x^k)_{k\in\NN}$ and
$\left(x^{\ell_k}\right)_{k\in\NN}$ be a convergent subsequence to
$\bar{x}$. Since $x^{k+1}$ belongs to $H(x^k,U^k)$, for all $k$, we
have
\begin{equation}\label{ptos-in-s*}
F(x^{\ell_k})+U^{\ell_k}\left(x^{\ell_k+1}-x^{\ell_k}\right)\dsty
\preceq_{K(F(x^{\ell_k}))} 0.\end{equation} By Assumption (A1),
Remark \ref{nota1} implies that $\left(U^{\ell_k}\right)_{k\in\NN}$
is bounded. So, the sequence
$\left(U^{\ell_k}(x^{\ell_{k}+1}-x^{\ell_k})\right)_{k\in\NN}$
converges to zero. By taking limits in \eqref{ptos-in-s*} and
recalling that $K$ is closed application, we obtain that
\begin{equation*}\label{escalarization}
\lim_{k\rightarrow\infty}F(x^{\ell_k})+U^{\ell_k}\left(x^{\ell_k+1}-x^{\ell_k}\right)=F(\bar{x})\in
-K(F(\bar{x})),
\end{equation*} implying that $\bar{x}\in S_*$\qed

Finally, we are ready to prove the convergence of the sequence
$(x^k)_{k\in\NN}$ generated by Algorithm $S$ to the solution which
lies closest to $x^0$.

\begin{theorem}  Define $x^*=P_{S_*}(x^0)$. Then $(x^k)_{k\in\NN}$ converges to $x^*$.
\end{theorem}
\proof   By  Lemma \ref{l:limitacao}, $(x^k)_{k\in\NN}\subset
B\left[\frac{1}{2}(x^0+x^*),\frac{1}{2}\rho\right]$ is bounded. Let
$\left(x^{\ell_k}\right)_{k\in\NN}$ be a convergent subsequence of
$(x^k)_{k\in\NN}$, and let $\bar{x}$ be its  limit. Evidently
$\bar{x}\in B\left[\frac{1}{2}(x^0+x^*),\frac{1}{2}\rho\right]$.
Furthermore, by Lemma \ref{l:optimalidad}, $\bar{x}\in S_*$. Since
$$S_* \cap  B
\left[\frac{1}{2}(x^0+x^*),\frac{1}{2}\rho\right]=\{x^*\},$$ and
recalling that $S_*$ is a convex and closed set,  we conclude that
$x^*$ is the unique limit point of $(x^k)_{k\in\NN}$. Thus,
$\left(x^{\ell_k}\right)_{k\in\NN}$ converges to $x^*\in S_*$. \qed

\section{Final Remarks}\label{sec-6}
In this paper we have presented two algorithms for finding a
solution to the $K$-convex variable inequalities problem. Using the
same hypotheses their convergence is shown. At Algorithm $S$ the
projection step involves more calculations than Algorithm $R$.
However, the sequence generated by the first algorithm has better
properties. In fact it converges to a solution of the problem, which
lies closest to the starting point. We emphasize that this last
special feature is  interesting  and it is useful in  specific
applications such as image reconstruction \cite{ImageR-1,ImageR-3}.
The main drawback of extending these algorithms to the infinite
dimensional spaces is that the existence of the subgradient has not
been shown in the variable order case.

\noindent As studied in \cite{gabi2,Gabriele-book}, variable orders can be
considered in two different ways,
$$y\preceq_K^1 \bar{y}\text{ if and only if }\bar{y}-y\in K(y)$$
or
$$y\preceq_K^2 \bar{y}\text{ if and only if }\bar{y}-y\in K(\bar{y}).$$
Problem \eqref{prob-1} corresponds with the inequalities defined by
$\preceq^1_K$. If the order is given by  $\preceq^2_K$, the
inequalities problem becomes
\begin{equation*}\text{find }x\in C\text{ such that }
F(x) \preceq_{K(0)} 0.
\end{equation*}
Since the cone $K(0)$ is fixed, the previous model is a non-variable
$K$-inequalities problem and it can be solved by the solution
algorithm proposed in  \cite{rob, yun-luis}.

We hope that this study will be useful for future research on other
more efficient variants of the subgradient iteration.
\subsubsection*{ACKNOWLEDGMENT}

This work started during a research stay of the first two authors at
National Institute for Pure and Applied Mathematics (IMPA) and it
was completed while the first author was visiting the University of
British Columbia. The authors are very grateful for the warm
hospitality of both institutions. The authors would like to thanks
to anonymous referees whose suggestions helped us to improve the
presentation of this paper.

%%%%%%%%%%%%%%%%%%%%%%%%
\end{document}